\documentclass[11pt]{amsart}
\usepackage{amsfonts}
\usepackage{amsmath}
\usepackage{amssymb}
\usepackage{amscd}
\usepackage{color}
\usepackage{amsthm}
\usepackage{indentfirst}
\usepackage[hmargin=3cm,vmargin=3cm]{geometry}

\newtheorem{theorem}{Theorem}

\newtheorem{corollary}{Corollary}[section]
\newtheorem{remark}[corollary]{Remark}
\newtheorem{definition}{Definition}

\newtheorem*{proposition*}{Proposition}
\newtheorem{thm}{Theorem}
\newtheorem*{thm*}{Theorem}
\newtheorem{prop}[corollary]{Proposition}
\newtheorem{lma}{Lemma}

\theoremstyle{definition}

\theoremstyle{remark}
\newtheorem*{pf}{Proof}

\newcommand{\R}{{\mathbb{R}}}
\newcommand{\Z}{{\mathbb{Z}}}
\newcommand{\C}{{\mathbb{C}}}

\newcommand{\D}{{\mathbb{D}}}

\newcommand{\del}{\partial}

\newcommand{\G}{\mathcal{G}}

\newcommand{\til}[1]{\widetilde{#1}}

\newcommand{\arr}[1]{\overrightarrow{#1}}

\newcommand{\vol}{{\rm vol}}

\newcommand{\om}{\omega}
\newcommand{\Om}{\Omega}

\newcommand{\eps}{\epsilon}

\newcommand{\cP}{\mathcal{P}}

\newcommand{\Ham}{{\rm Ham}}

\begin{document}

%
%
%
%
%
%
%
%

\title{Proof of the main conjecture on $g$-areas}

\author{Fran\c{c}ois Lalonde}
\address[Fran\c{c}ois Lalonde]{D\'epartement de math\'ematiques et de Statistique, Universit\'e de Montr\'eal, C.P. 6128, Succ. Centre-ville, Montr\'eal H3C 3J7, Qu\'ebec, Canada}
\email{lalonde@dms.umontreal.ca}

\author{Egor Shelukhin}
\address[Egor Shelukhin]{Centre de Recherches Math\'ematiques, Universit\'e de Montr\'eal, C.P. 6128, Succ. Centre-ville, Montr\'eal H3C 3J7, Qu\'ebec, Canada}
\email{egorshel@gmail.com}



\begin{abstract}

{\em In this paper, we prove the main conjecture on $g$-areas arising from \cite{LT}. That conjecture was announced by the first author in 2004. It states that the $g$-area of any Hamiltonian diffeomorphism $\phi$ is equal to the positive Hofer distance between $\phi$ and the subspace of Hamiltonian diffeomorphisms that can be expressed as a product of at most $g$ commutators. 
}

\end{abstract}

\maketitle

   \section{Introduction}
   
      An interesting direction in Symplectic Topology is the study, on the infinitesimal level, of the lack of integrability and, on the global level, of the lack of commutativity of the main objects considered. One of the problems in this direction is the ``commutator problem for Hamiltonian diffeomorphisms'',  wherein one is interested in understanding to which extent the subspace of the space of Hamiltonian diffeomorphisms that consists of all products of $k$ commutators catches the topology and the geometry of the full space of Hamiltonian diffeomorphisms in the Hofer norm. Note that we work here in the Hofer metric because it is essentially the only one which is intrinsically defined on the group of Hamiltonian diffeomorphisms.
      
     So let $(M, \om)$ be a closed symplectic manifold, that is to say a real smooth manifold (i.e. of class $C^{\infty}$) of real dimension $2n$ endowed with a symplectic form: a real smooth differential 2-form with values in the real numbers which is non-degenerate on the real tangent space at any point $p \in M$, and whose exterior derivative is closed. Define the group of Hamiltonian diffeomorphisms as the set of diffeomorphisms of $M$ that can be expressed as the time-one map of the Hamiltonian evolution equation of a real smooth time-dependent function
     
     $$
     H: M \times [0,1]  \to \R.
     $$
     
    More precisely, each smooth function $H: M \to \R$ gives rise to a symplectic gradient, i.e a vector field $X_H$ on $M$ defined by using the non-degenerate 2-form $\omega$ to pair $-dH$ with a vector (exactly as in the case of the scalar product, which in this latter case leads to the usual Riemannian gradient). These functions represent the ``total energy'' and therefore bare usually the letter $H$ for Huyghens who discovered the concept of total energy of a physical system (this terminology was introduced in Lagrange's book of Analytical Mechanics). Now, if $H$ depends on the $t$-variable as above, this procedure produces a time-dependent vector field. Integrating that time-dependent ODE from time $0$ to $1$ gives a diffeomorphism. It is not difficult to see that it preserves the symplectic structure. These symplectic diffeomorphisms are called the {\it Hamiltonian diffeomorphisms} of the manifold $M$.
     
     By a deep theorem of Banyaga \cite{B}, the group of Hamiltonian diffeomorphisms $\Ham(M)$ of any closed symplectic manifold $M$ is simple. So it does not admit any intrinsically defined (and hence normal) non-trivial subgroup. Therefore the subset of all Hamiltonian diffeomorphisms $\phi$ that can be expressed as a finite product of commutators of elements in the group of Hamiltonian diffeomorphisms
     $$
     \phi = [f_1, g_1] \ldots [f_k, g_k] \quad \quad f_i, g_i \in \Ham(M)
     $$
     
     \noindent
     (where $k$ varies) 
     
     \medskip
     \noindent
     must, as it is non-trivial, be the full group of Hamiltonian diffeomorphisms. Here the commutator $[f,g]$ is simply 
     
     $$f \circ g \circ f^{-1} \circ g^{-1}.$$ 
     
     Thus:
     
\begin{proposition*} Each Hamiltonian diffeomorphism of a closed symplectic manifold can be expressed as a finite product of commutators of elements in the group of Hamiltonian diffeomorphisms.
\end{proposition*}

   In this article, we explore the relation between three objects: (1)  the subspace of $\Ham(M)$ consisting of all Hamiltonian diffeomorphisms that can be expressed as a product of at most $k$ commutators, that we will denote by $C_k(M)$;  (2) the positive Hofer norm on $\Ham(M)$; and (3) the $g$-area introduced in \cite{LT}, where $g$ stands for {\it genus}.
   
   Let's first define the positive Hofer norm. It associates to each Hamiltonian diffeomorphism $\phi$ of a closed symplectic manifold a real number in $[0, \infty)$ that we will denote by $\| \phi \|_H^+$. It is an infimum over a set of geometric objects that we now define (see \cite{LM} for more details): 
   
   \begin{definition} A Hamiltonian fibration is by definition a fiber bundle 
   
   $$\pi : P \to B$$
   where $B$ is a smooth compact manifold with or without boundary, $P$ is a symplectic manifold with symplectic structure $\Om$, and $\pi$ is a smooth map such that the restriction of $\Om$ to every fiber $P_b:=\pi^{-1}(b \in B)$  is non-degenerate (and therefore induces a symplectic structure on each fiber). We denote by $(M, \omega)$ a fixed symplectic manifold symplectomorphic to the fiber endowed with the restriction of $\Om$ (which is well-defined, since all fibers are symplectomorphic). 
   
   \end{definition}
   
   \begin{remark}
Let $\gamma$ be a loop in $B$ based at $b \in B.$ Fixing a symplectomorphism \[{I_{b}:(P_b,\Om|_{P_b}) \to (M,\omega)}\] we define the monodromy of $(P,\Om)$ over $\gamma$ to be symplectomorphism \[Hol_{(P,\Om)}(\gamma):(M,\omega) \to (M,\omega),\] defined by \[Hol_{(P,\Om)}(\gamma) = I_{b} \Gamma_\gamma (I_{b})^{-1},\] where $\Gamma_\gamma$ is the parallel transport along $\gamma$ of the Ehresmann connection on $\pi : P \to B$ given by the distribution $H_p:=\mathrm{Ker} (\Om_p)$ for each $p \in P.$ In what follows, when fixing a point $b \in B$ we always fix an identification $I_{b}$ as above and consider it part of the data.
   \end{remark}
   
      Note that if the base $B$ is a real oriented surface with boundary consisting of a finite union of loops, then the restriction of $\Om$ to each codimension $1$ submanifold of the form $\pi^{-1}(\gamma)$, where $\gamma$ is any connected component of the boundary of the base $B$, has a real dimension $1$ kernel lying on $\pi^{-1}(\gamma)$. That kernel induces a flow on the boundary whose monodromy around $\gamma$ is a Hamiltonian diffeomorphism of the fiber.
   
     Then we have:
     
  \begin{definition} Let $\phi$ be a Hamiltonian diffeomorphism of a closed symplectic manifold $(M, \om)$. Consider the Hamiltonian fibrations $P \to D^2$ over the closed unit $2$-disc, with fiber $(M, \om)$, such that  the monodromy round the boundary of the disc based at a fixed base-point is equal to $\phi$. Then define the positive Hofer norm $\| \phi \|_H^+$ as the infimum of the quantity
  
  $$
     A(P,\Om) = \frac{\vol(P, \Om)}{\vol(M,\om)}
     $$
     over all such such fibrations.  
      \end{definition}

     We can extend this to a pseudo-distance by setting $d_0(\phi, \psi) = \| \phi \circ \psi^{-1} \|_H^{+}$.

     \begin{remark}
     We denote by $\cP(D,\phi)$ the class of such Hamiltonian fibrations. We note that choosing trivializations of $P \to D$ standard over the base-point, we can define the class $\cP(D,\til{\phi})$ where the holonomy path with respect to any trivialization over the boundary lies in class $\til{\phi} \in \til{\G}.$ (Note that the class in $\til{\G}$ of the holonomy path does not depend on the trivialization.) Then for $\til{\phi} \in \til{\G}$ we define \[\| \til{\phi}  \|_H^{+} = \displaystyle{\inf \{A(P,\Om)\,| \;{(P,\Om) \in \cP(D,\til{\phi})}} \},\] and \[d_0(\til{\phi},\til{\psi}) = \| \til{\phi} \til{\psi}^{-1} \|_H^{+}.\]
     \end{remark}
      
    \begin{remark}
   See \cite{M} for an alternative definition of the positive Hofer norm and for the equivalence of the two definitions.
    \end{remark}  
    
    Now let's consider a different ``norm" that we will call the $g$-{\it area}, introduced in \cite{LT}. It is defined in the following way.
    
    \begin{definition} Let $\phi$ be a Hamiltonian diffeomorphism of a closed symplectic manifold $(M, \om)$. Consider the Hamiltonian fibrations $P \to \Sigma_g$, with fiber $(M, \om)$ over the compact orientable real surface of genus $g$ with one boundary component equal to a circle. Thus $\Sigma_g$ is obtained from the closed surface of genus $g$ by removing one open disc. Equivalently, it is the surface obtained by adding $g$ handles to the closed unit disc. Assume that  the monodromy round the boundary of $\Sigma_g$ is equal to $\phi$. Then define the $g$-area of $\phi$ as the infimum of the quantity
  
  $$
     A(P,\Om) = \frac{\vol(P, \Om)}{\vol(M,\om)}
     $$
     
      over all such such fibrations.  
      \end{definition}
      
      \begin{remark}
	We denote by $\cP(\Sigma_g,\phi)$ the class of such Hamiltonian fibrations.
      \end{remark}
      
      Our main theorem is:
      
      \begin{theorem}\label{Theorem: main} Let $\phi$ be a Hamiltonian diffeomorphism of a closed symplectic manifold $(M, \om)$. Then the $g$-area of $\phi$ is equal to the infimum of the positive Hofer distance between $\phi$ and $\psi$, the infimum being taken over all $\psi$ in the commutator subspace $C_g(M)$.
      \end{theorem}
      
        Thus the genus $g$-area measures exactly the usual positive Hofer distance to the commutator subspace $C_g(M)$. Half of this result was obtained in \cite{LT}: it was proved there that the genus $g$-area is bounded above by the positive Hofer distance between $\phi$ and $C_g(M),$ sharpening a previous result of Entov \cite{E}. The result of \cite{LT} implies in particular that if a diffeomorphism belongs to $C_k(M)$, then its genus $g$-areas must vanish for all $g \le k.$

\section{Structure of the proof.}

For brevity, put $\G = \Ham(M)$ for the group of Hamiltonian diffeomorphisms of a closed symplectic manifold $(M,\om),$ and for an integer $k \geq 0$ let $C_k = C_k(M) \subset \G$ be the subset of all Hamiltonian diffeomorphisms that can be written as a product of at most $k$ commutators. For an element $c \in \G$ denote by $\Lambda_c$ its conjugacy class in $\G.$ For any pseudo-distance $d_*$ on $\G,$ and $\phi \in \G$ put $d_*(\phi,C_k) := \inf_{\psi \in C_k} d_*(\phi,\psi)$ for the $d_*$-distance between $\phi$ and $C_k.$


   Before beginning the proof, let's generalize the definition of our genus $g$-areas in the following way.
   
\begin{definition}\label{Definition: general} Let $M$ be a symplectic manifold, and let  $\phi_1, \ldots, \phi_k$ be Hamiltonian diffeomorphisms of $M$. Let $g \ge 0$ be a given positive integer. We denote by $\Sigma_{g,k}$ the real oriented surface of genus $g$ with $k$ marked ordered ends, that is to say the real closed oriented surface of genus $g$ with $k$ ordered open discs removed and a choice of a point $z_j$ on its boundary, for each $j= 1, \ldots,k$ . We will denote that simply by  $\Sigma$ if there is no confusion. We denote by $\cP(\Sigma_{g,k},\phi_1, \ldots, \phi_k)$ the class of Hamiltonian fibrations $M \subset P \to \Sigma_{g,k}$  with holonomy $\phi_j$ over the $j^{th}$ end for all $j= 1, \ldots, k$ (the holonomy is taken with respect to the base point $z_j$).

    Then define the genus $g$-area of the sequence $\phi_1, \ldots, \phi_k$ by taking the infimun of 
    
      $$
     A(P,\Om) = \frac{\vol(P, \Om)}{\vol(M,\om)}
     $$
over all such Hamiltonian fibrations. 

\end{definition}


Now, the proof of Theorem \ref{Theorem: main} rests on an additional pseudo-metric on $\G$ and $\til{\G}$ defined by taking the infimum of areas for Hamiltonian fibrations over the annulus with fixed holonomies on the boundary components. Put $A=[0,1] \times S^1$ for the annulus. The following special case of Definition \ref{Definition: general} is central to our proof.

\begin{definition}
For $\phi,\psi \in \G$ define the class $\cP(A,\phi,\psi)$ of Hamiltonian fibrations over $A$ with holonomy $\phi$ over the boundary component $S_0 = \{0\} \times S^1$ and holonomy $\psi$ over $S_1 = \{1\} \times S^1.$ Similarly, choosing trivializations fixing the base-points on the boundary components, we have the class $\cP(A,\til{\phi},\til{\psi}).$ 

Then put \[d_A(\phi,\psi)= \displaystyle{\inf \{A(P,\Om)\,| \;{(P,\Om) \in \cP(A,\phi,\psi)}}\}, \] and 
\[d_A(\til{\phi},\til{\psi})= \displaystyle{\inf \{A(P,\Om)\,| \;{(P,\Om) \in \cP(A,\til{\phi},\til{\psi})}}\}.\]
\end{definition}

\begin{remark}
We note that in the definition of $\cP(A,\til{\phi},\til{\psi})$ the classes in $\til{\G}$ of the holonomies of the boundary components do depend on the trivialization. However, we can either work with an a-priori fixed trivialization varying only the connection and compatible symplectic form, or note that since the area $A(P,\Om)$ is invariant under gauge transformations and changing the trivialization changes the holonomies over the boundary components by the same homotopy class of loops, we get the same value of $d_A(\til{\phi},\til{\psi})$ for the two classes of Hamiltonian fibrations defined by any two different trivializations. Indeed we can always correct the homotopy classes of the holonomy paths in $\til{\G}$ by a gauge transformation (see proof of Lemma \ref{Lemma: contractible existence}).
\end{remark}

We have the following invariance properties of $d_A$ and relations between $d_A$ and $d_0.$ 

\begin{prop}\label{Proposition: d_A and d_0 from 1}
For all $\phi \in \G$ we have \[d_A(\phi,1) = d_0(\phi,1).\] An analogous statement holds for $\til{\G}.$ Namely \[d_A(\til{\phi},\til{1}) = d_0(\til{\phi},\til{1})\] for all $\til{\phi} \in \til{\G},$ where $\til{1} \in \til{\G}$ is the class of contractible loops.
\end{prop}

\begin{prop}\label{Proposition: invariance of d_A}
For all $\phi,\psi,h \in \G,$ \[d_A(\phi,\psi \cdot \Lambda_h) = d_A(\phi \cdot \Lambda_{h^{-1}},\psi).\] An analogous statement holds for $\til{\G}.$ Namely \[ d_A(\til{\phi}, \til{\psi} \cdot \Lambda_{\til{h}}) = d_A(\til{\phi} \cdot \Lambda_{\til{h}^{-1}}, \til{\psi}) \] for all $\til{\phi},\til{\psi},\til{h} \in \til{\G}.$
\end{prop}

Combining Propositions \ref{Proposition: d_A and d_0 from 1} and \ref{Proposition: invariance of d_A} we have 
\begin{prop}\label{Proposition: d_A and d_0 for congugacy classes}
For all $\phi,\psi \in \G$ we have $d_A(\phi,\Lambda_\psi) = d_0(\phi,\Lambda_\psi).$
\end{prop}

Indeed \[d_A(\phi,\Lambda_\psi) = d_A(\phi\cdot \Lambda_{\psi^{-1}}, 1) = d_0(\phi\cdot \Lambda_{\psi^{-1}}, 1) = d_0(\phi,\Lambda_\psi)\] for all $\phi,\psi \in \G.$

In \cite{LT} the following inequality was shown.

\begin{thm*}\label{Theorem: upper bound on g-area}\cite{LT}
For all $\phi \in \G$ \[|\phi|_g \leq d_0(\phi,C_g).\]
\end{thm*}

Here we prove the reverse inequality

\begin{prop}\label{Prop: upper bound on g-area}
\[|\phi|_g \geq d_0(\phi,C_g),\] for all $\phi \in \G.$
\end{prop} 

This establishes the main result of this note:

\begin{thm}\label{Theorem: equality on g-area}
For all $\phi \in \G$ \[|\phi|_g = d_0(\phi,C_g).\]
\end{thm}

One proves Proposition \ref{Prop: upper bound on g-area} as follows.

\begin{pf}
First, by Proposition \ref{Proposition: d_A and d_0 for congugacy classes} and the fact that $C_g$ is conjugation invariant, it immediately follows that it is enough to prove that for all $\phi \in \G,$

\begin{equation}\label{Equation: actual lower bound}
|\phi|_g \geq d_A(\phi,C_g).
\end{equation}

To do this we note that given a Hamiltonian fibration $(P,\Om) \in \cP(\Sigma_g, \phi)$ of area $A(P,\Om)$ we can construct a Hamiltonian fibration $(P',\Om') \in \cP(A, \phi, \psi)$ for $\psi \in C_g$ by simply cutting $\Sigma_g$ along the $2g$ loops $\{\alpha_j,\beta_j\}_{1 \leq j \leq g}$ (geodesic loops with respect to the hyperbolic metric) and obtaining a punctured fundamental domain in the upper half-plane, isomorphic to $A,$ from the punctured surface of genus $g,$ isomorphic to $\Sigma_g.$ The new fibration is formally the pull-back $(P',\Om') = \pi^*(P,\Om)$ of the original fibration with respect to the natural surjective map $\pi$ from the fundamental domain to the surface. Then denoting $a_k := Hol_{(P,\Om)}(\alpha_k),\; b_k := Hol_{(P,\Om)}(\beta_k)$ for all $1 \leq k \leq g$ we see that the holonomy of $(P',\Om')$ along the new boundary component is $\psi = [a_1,b_1]\cdot ... \cdot [a_g,b_g] \in C_g.$ Moreover $A(P',\Om') = A(P,\Om).$ Taking infima then finishes the proof.

\end{pf}

\section{Proofs}\label{Section: Proofs}

In this section we prove Propositions \ref{Proposition: d_A and d_0 from 1} and \ref{Proposition: invariance of d_A}. Both are proved using the following general technique for gluing Hamiltonian fibrations that was first introduced in \cite{LMglueing} and later developed in \cite{L}.

\begin{lma}\label{Lemma: Gluing along same holonomy} 
Let $\phi_1, \ldots, \phi_k, \phi_1', \ldots, \phi_{k'}'$ be Hamiltonian diffeomorphisms of the same symplectic manifold $(M, \om)$. Let  $(P,\Om) \in \cP(\Sigma_{g,k},\phi_1, \ldots, \phi_k)$ and $(P',\Om') \in \cP(\Sigma_{g',k'},\phi_1', \ldots, \phi_{k'}')$ be two Hamiltonian fibrations with same fiber $(M, \om)$, and suppose that some monodromy in the first fibration is equal to some monodromy in the second one -- we may assume, after reordering that $\phi_1 = \phi_1'$. Then, for any $\delta > 0$, one can glue the two fibrations along their common monodromy and obtain a fibration $(P'',\Om'') \in \cP(\Sigma_{g+g',k+k'-2},\phi_2, \ldots, \phi_k, \phi_2', \ldots, \phi_{k'}')$ of area $A(P'',\Om'')$ less than or equal to $A(P,\Om) + A(P',\Om') + \delta$.
\end{lma}

\begin{pf}

   Denote by $\pi, \pi'$ the two projections of the Hamiltonian fibrations.  Let $End_1$ denote the end of the closed oriented surface where the monodromy $\phi_1$ occurs, and $C_1$ its boundary. Let $R_1$ denote the restriction of the fibration $P$ to its boundary component $C_1$. Define similarly, with prime symbols, the same objects on the second fibration $P'$.   We first define an identification of  $R_1$ with $R_1'$ in such a way as to respect the fibers symplectically and the characteristic foliations of $\Om$ on $R_1$ and $\Om'$ on $R_1'$. This is easy if we transform the characteristic foliations into Hamiltonian vector fields whose flows respect the fibers. For this,  take a real function $F$ defined on $End_1$ in such a way that, in some given trivialization of $End_1$ $End_1 \to (-\eps, 0] \times S^1$, the function  $F$ is the projection $(-\eps, 0] \times S^1 \to (-\eps, 0] $. Then smooth out this function near $- \eps$ so that it coincide with a constant in some small open neighbourhood  of $- \eps$ inside $(-\eps, 0] $. Define $H$ as the pull back to the fibration of that function, i.e $H = F \circ \pi$, and denote by $X_H$ its symplectic gradient. Its flow $\{\psi_t\}_{t \in [0,1]}$ is Hamiltonian and therefore symplectic and it preserves the fibers of the fibration since it comes from the pull-back of a function in the base. Do the same for the second fibration. 
    
     Now pick a base point $b$ on the boundary component $C_1$ in the base corresponding to the first end $End_1$  of the first fibration and pick similarly a point $b'$ on the boundary component $C_1'$ in the base corresponding to the first end $End_1'$ of the second fibration. Now identify the symplectic fiber $P_b$ of $P$ over $b$ with the symplectic fiber $P_{b'}'$ of $P'$ over $b'$ using the fixed identifications of these fibers with $(M,\om),$ and call the resulting symplectic diffeomorphism $I$. Extend that symplectic identification to a symplectic diffeomorphism $\Psi: R_1 \to  R_1'$ by matching the characteristic flow $X_H$ with $X_{H'}$. More precisely, both Hamiltonian diffeotopies, $\psi_t$ and $\psi'_t$ , close up in time one because they come from $H, H'$ which both come from standard functions $F, F'$ defined on the annuli on the base (and both annuli are parametrized similarly). They close up to the same monodromy $\phi$. So we define $\Psi$ as the map that sends the Hamiltonian flow of $H$ to the Hamiltonian flow of $H'.$ That is to say, we set, for $t$ on $C_1$ and $x$ in the fiber over $t$, and $I$ the identification of the fibers at the base points
$$
\Psi(x,t) =  \psi'_t(x,t)  \circ I \circ \psi_t^{-1} (x,t) 
$$ 
which is symplectic. Obviously, it is well define since, after one complete turn, the maps $\psi_1$ and $\psi'_1$ (i.e the monodromies) coincide.
       
     By the normal form of symplectic structures near coisotropic submanifolds, $\Psi$ can be extended to a symplectic diffeomorphism
     $$
     \Phi: R_1^+  \to (R_1')^-
     $$
     where the source here is an open outward collar neighbourhood of $R_1$ and the target is an open inward collar neighbourhood of $R_1'$ (actually, this diffeomorphisms is provided by the flows of $H, H'$ themselves since they are defined on such collar neighbourhoods). 

    This leads, by glueing using $\Phi$, to a new total space, and we just have to show that this can be made into a Hamiltonian fibration. The point here is to get a fibration with non-degenerate fibers. So assume that on the glued part, we take as fibers those coming from $R_1^+$. Transporting them through $\Phi$ gives non-degenerate fibers in $(R_1')^-$ that coincide with the fibers of $P'$ over $C_1'$. But if the construction of $\Phi$ is realized on a sufficiently small neighbourhood, the image by $\Phi$ of each fiber over $(s,t)$ is sufficiently $C^1$-close to the fiber of $(R_1')^-$ at $(-s,t)$ so that it can be viewed as the graph of a $C^1$-small function to the real numbers (keep in mind that this push-forward of the fibers respect time). Then isotop that graph to the graph of the zero function over a small inner neighbourhood, so that these fibers coincide over $\{-\eps'\} \times S^1$ with the fibers of $P'$. 
    
   It is obvious, by construction, that the area of the new fibration obeys the upper bounds stated in the theorem.

\end{pf}

\smallskip
\noindent
{\em Remark.}  Note that the lemma is obviously true if one glues $m$ ends of the first fibration with $m$ ends of the second one, assuming of course that their monodromies match pairwise. This would lead to a Hamiltonian fibration
$P'' \in \cP(\Sigma_{g+g' + m-1,k+k'-2m},\phi_{m+1}, \ldots, \phi_k, \phi_{m+1}', \ldots, \phi_{k'}')$ with area satisfying the same upper bound.

\medskip
We will also be using the following two lemmas on the existence of Hamiltonian fibrations of arbirtrarily small areas with prescribed holonomy paths over the boundary.

\begin{lma}\label{Lemma: contractible existence}
Given a contractible loop $\arr{\phi} = \{\phi_t\}$ in $\G,$ for every $\eps > 0$ there exists a Hamiltonian fibration $(P^\#,\Om^\#)$ of area $A(P^\#,\Om^\#) < \eps$ over the disc $D$ with the holonomy path $\arr{\phi}$ over the boundary $S^1.$
\end{lma}

\begin{lma}\label{Lemma: square existence}
Given two paths $\arr{e}=\{e_t\}$ and $\arr{f}=\{f_t\}$ in $\G$ starting at $1,$ for every $\eps > 0$ there exists a Hamiltonian fibration $(P^\#,\Om^\#)$ of area $A(P^\#,\Om^\#) < \eps$ over $[0,1]\times [0,1]$ with holonomy paths $\arr{e}, e_1 \arr{f} e_1^{-1}, \arr{e}, \arr{f}$ over the boundary paths $\gamma_0=\{(0,t)\}_{t \in [0,1]},\; \eta_1=\{(t,1)\}_{t \in [0,1]},\; \gamma_1=\{(1,t)\}_{t \in [0,1]},\; \eta_0=\{(t,0)\}_{t \in [0,1]}$ respectively.
\end{lma}

We defer the proofs of these two lemmas to the end of the section and proceed to prove Propositions \ref{Proposition: d_A and d_0 from 1} and \ref{Proposition: invariance of d_A}.
\medskip

\begin{pf}\emph{(Proposition \ref{Proposition: d_A and d_0 from 1})}

Assume that we know the analogue of this statement for $\til{\G}.$ Namely \[d_A(\til{\phi},\til{1}) = d_0(\til{\phi},\til{1})\] for all $\til{\phi} \in \til{\G},$ where $\til{1} \in \til{\G}$ is the class of contractible loops. Then we note that \[d_A(\phi,1) = \inf_{\til{\phi} \in \pi^{-1}(\phi),\gamma \in \pi_1(\G)} d_A(\til{\phi},\gamma),\] and a similar formula holds for $d_0(\phi,1).$ However since $\pi_1(\G) \subset \mathrm{Z}(\til{\G}),$ we have $\Lambda_{\gamma} = \{\gamma\}$ and therefore by the analogue of Proposition \ref{Proposition: invariance of d_A} for $\til{\G},$ \[d_A(\til{\phi},\gamma) = d_A(\til{\phi}\gamma^{-1},\til{1}).\] Therefore \[d_A(\phi,1) = \inf_{\til{\phi} \in \pi^{-1}(\phi)} d_A(\til{\phi},\til{1}) = \inf_{\til{\phi} \in \pi^{-1}(\phi)} d_0(\til{\phi},\til{1}) = d_0(\phi,1),\] by the bi-invariance of the positive Hofer metric.

The inequality $d_A(\til{\phi},\til{1}) \leq d_0(\til{\phi},\til{1})$ is easy. Simply note that given a fibration in $\cP(D,\til{\phi}),$ one can trivialize it over a small neighbourhood of the base-point and cut out its restriction to a small open disc in this neighbourhood. This gives us a fibration $\cP(A,\til{\phi},\til{1}),$ with smaller area. 

We prove the reverse inequality $d_A(\phi,1) \leq d_0(\phi,1).$ Take a fibration  $(P,\Om) \in \cP(A,\til{\phi},\til{1}).$ Let $c$ denote its contractible holonomy loop over $S_1.$ Given an $\epsilon >0,$ Lemma \ref{Lemma: contractible existence} gives us a fibration $(P^\#,\Om^\#)$ over a disc with $A(P^\#,\Om^\#) < \epsilon/2$ and the boundary holonomy path $c.$ Gluing $(P,\Om)$ and $(P^\#,\Om^\#)$ along the corresponding boundary components, given by Lemma \ref{Lemma: Gluing along same holonomy}, gives us a bundle $(P',\Om') \in \cP(D,\til{\phi})$ with $A(P',\Om') < A(P,\Om) + A(P^\#,\Om^\#) + \epsilon/2 < A(P,\Om) + \epsilon.$ Taking infima now finishes the proof.
\end{pf}

\begin{pf}\emph{(Proposition \ref{Proposition: invariance of d_A})}
Consider a Hamiltonian fibration $(P,\Om) \in \cP(A,\phi,\psi h')$ for $h' \in \Lambda_h.$ Consider the path $\gamma = \{(0,t)\}_{t \in [0,1]}$ between the marked points $(0,0)$ and $(0,1)$ on the boundary components $S_0 = S^1 \times \{0\}$ and  $S_1 = S^1 \times \{1\}$ of  $A = S^1 \times [0,1].$ Denote by $c:[0,1] \to \G$ the holonomy path of $(P,\Om)$ over $\gamma$ (with respect to any trivialization of $P$).  Cut $A$ along $\gamma$ to get a Hamiltonian fibration $(P^0,\Om^0)= \pi^* (P,\Om)$ over $[0,1] \times [0,1]$ for the natural surjection $[0,1] \times [0,1] \to A.$
Its holonomy paths over the vertical boundary paths $\gamma^0_0$ and $\gamma^0_1$ are equal to $c=\{c_t\}.$ Consider the Hamiltonian fibration $(P^\#,\Om^\#)$ over $[0,1] \times [0,1]$ provided by Lemma \ref{Lemma: square existence} with $e=c$ and $f=\{f_t\}$ with $f_t = (h'_t)^{-1}$ for a path $\{h'_t\}$ with $h'_0 = 1$ and $h'_1 =  h'.$ Its holonomies over $\gamma_0$ and $\gamma_1$ are equal to $c.$ Hence we can glue $(P_0,\Om_0)$ and $(P^\#,\Om^\#)$ along the pairs of paths $\gamma^0_1, \gamma_0$ and $\gamma_1,\gamma^0_0,$ and obtain a bundle $(P',\Om') \in \cP(A,\phi (h'')^{-1}, \psi),$ where $(h'')^{-1} = c_1^{-1} (h')^{-1} c_1 \in \Lambda_{h^{-1}},$ over $A$ with $A(P',\Om') < A(P,\Om) + \eps.$ Taking infima proves $d_A(\phi \cdot \Lambda_{h^{-1}}, \psi) \leq d_A(\phi, \psi \cdot \Lambda_{h}).$
The other direction is similar.

This proof clearly goes through for the analogous statement for $\til{\G}.$ 
\end{pf} 

\begin{pf}\emph{(Lemma \ref{Lemma: contractible existence})}
For concreteness, take the two-disc $D$ to be $\D=\{|z|\leq 1\} \subset \C.$ Consider the inverse contractible loop $\arr{\psi}=\{\psi_t\}_{t \in S^1},$ with $\psi_t = (\phi_t)^{-1}$ for all $t \in S^1 = \R/\Z,$ and the map $\Psi:\D \times M \to \D \times M$ given by $\Psi(z,x)=(z,\psi_{z}\cdot x)$ where the map $\D \to \G$ given by $z \mapsto \psi_z$ is a contraction of $\psi,$ namely when rewritten in polar coordinates $(s,t)$ with $z=s \cdot e^{(i 2\pi) t},$ we have $\psi_t = \psi_{1,t}$ for all $t \in S^1.$

Endow $\D \times M$ with the structure of a trivial Hamiltonian fibration $P^{triv} \to \D,$ with the trivial connection. Consider the pull-back Hamiltonian fibration with connection $P^\#=\Psi^*(P^{triv})$ with respect to $\Psi.$ We claim that its holonomy path over $S^1 \cong \del \D$ is $\arr{\phi}$ and that for every $\eps >0$ there is a compatible symplectic form $\Om^\#$ on $P^\#$ with $A(P^\#,\Om^\#) < \eps.$

For the holonomy we work in polar coordinates as above, and partially compute the differential of the map $\Phi$ near $\del (P^{triv}).$

\begin{align}
  D(\Phi)(s,t,x)\,\del_t &= \del_t + Y_t\circ \psi_t(x) \label{Equation: differential of Psi 1}\\
  D(\Phi)(s,t,x) \,\xi &= D(\psi_{s,t})(x) \,\xi , \;\;\; \xi \in T_x M. \label{Equation: differential of Psi 2}
\end{align}

Here $Y_t$ is the Hamiltonian vector field generating $\psi.$ We note that the horizontal vector field in $P^{triv}$ covering $\del_t$ on $\del D$ is $\del_t$ as a vector field on $\del(P^{triv}).$ Hence, the horizontal vector field in $P^\#$ covering $\del_t$ on $\del_D$ will be $D(\Psi^{-1})(\del_t).$ We compute by Equation \ref{Equation: differential of Psi 1}: \[(D(\Psi^{-1})\del_t)(1,t,x) = (D\Psi(1,t,x))^{-1} \del_t = \del_t - (D\Psi(1,t,x))^{-1} Y_t\circ \psi_t(x).\] 

And by Equation \ref{Equation: differential of Psi 2} we conclude \[(D(\Psi^{-1})\del_t)(1,t,x) = \del_t - (D(\psi_t)(x))^{-1} Y_t\circ \psi_t(x).\]

But $X_t = - (D(\psi_t)(x))^{-1} Y_t\circ \psi_t(x)$ is the vector field generating $\phi,$ hence the horizontal vector field covering $\del_t$ in $P^\#$ is simply \[\del_t + X_t,\] and hence the holonomy path of $P^\#$ over $\del \D \cong S^1$ is $\arr{\phi}.$

For the size we now choose a symplectic form $\Om_\eps = \eps \cdot \om_\D \oplus \om_M$ on $P^{triv}$ compatible with the trivial connection. Here $\om_\D$ is the standard symplectic form on $\D \subset \C$ and $\om_M$ is the chosen symplectic form on $M.$ Then $\Om^\#_\eps = \Psi^* \Om_\eps$ will be a symplectic form compatible with the pull-back connection on $P^\#.$ It is sufficient to show that $A(P^\#,\Om^\#_\eps) \to 0$ as $\eps \to 0.$ Writing $z = a+ib,$ we put $V_z = (\del_a \psi_z)\circ (\psi_z)^{-1}$ and $W_z = (\del_b \psi_z)\circ (\psi_z)^{-1}$ for the deformation vector fields in the $a$ and $b$ directions, and let $F_z$ and $G_z$ be the corresponding normalized generating Hamiltonians. Calculating similarly to Equations \ref{Equation: differential of Psi 1} and \ref{Equation: differential of Psi 2} we have \[\Om^\#_\eps = \Psi^* \Om_\eps = \eps da\wedge db + d(F_z \circ \psi_z) \wedge da + d(G_z \circ \psi_z) \wedge db + \om_M.\] Hence \begin{align*}(\Om^\#_\eps)^{n+1} &= \eps (n+1)(\om_M)^n da\wedge db - 2 {n+1 \choose 2}(\om_M)^{n-1} \cdot d(F_z \circ \psi_z) \wedge d(G_z \circ \psi_z) \wedge da \wedge db \\ &= 
(\eps (n+1) + \{F_z,G_z\} \circ \psi_z) (\om_M)^n da\wedge db.\end{align*} Therefore, since Poisson brackets integrate to zero against the symplectic volume form, we get \[A(P^\#,\Om^\#_\eps) = \eps \pi (n+1),\] which tends to zero as $\eps \to 0.$ This finishes the proof.
\end{pf}

\begin{pf}(Lemma \ref{Lemma: square existence})
It suffices to note that $\arr{e} \# (e_1 \arr{f} e_1^{-1}) e_1 \# (\arr{e}^{-1}) e_1 f_1 \# (\arr{f}^{-1}) f_1$ is a contractible loop based at $1,$ and use Lemma \ref{Lemma: contractible existence}. Here $\arr{e}^{-1} =\{e_t^{-1}\}_{t \in [0,1]}$ and $\arr{f}^{-1} =\{f_t^{-1}\}_{t \in [0,1]}.$
\end{pf}
\bibliographystyle{amsplain}

\begin{thebibliography}{999}

\bibitem{B} A. Banyaga, Sur la structure du groupe des diff\'eomorphismes qui p\'eservent une forme symplectique, {\it Commentarii Mathematici Helveticii} {\bf 53} (1978), 174 -- 227.

\bibitem{E} M. Entov, Commutator length of symplectomorphisms, {\it Commentarii Mathematici Helvetici} {\bf 79} (2004), 58 -- 104.

\bibitem{L} F. Lalonde, A field theory for symplectic fibrations over surfaces with applications, {\it Geometry and Topology}  {\bf 8} (2004), 1189 -- 1226.  

\bibitem{LMglueing}  F. Lalonde and D. McDuff, Hofer's $L^{\infty}$-geometry:
energy and stability of Hamiltonian flows part II, {\it Inventiones Mathematicae} {\bf 122} 
(1995), 35 -- 69.

\bibitem{LM} F. Lalonde and D. McDuff, Symplectic stuctures on fiber bundles, {\it Topology} {\bf 42} (2003), 309 -- 347.


\bibitem{LT} F. Lalonde and A. Teleman, g-areas and commutator length, { \it International Journal of Mathematics} {\bf 24} (2013).

\bibitem{M} D. McDuff, Geometric Variants of the Hofer norm, {\it Journal of Symplectic Geometry}, {\bf 1} (2002), 197--252

\end{thebibliography}

\end{document}